\documentclass[a4paper,10pt]{article}

\usepackage{amsfonts,amssymb,mathrsfs,amscd}
\def \R {{\mathbb {R}}}

\def \Z {{\mathbb {Z}}}

\def\eps{\varepsilon}

\title{\bf Recurrence of integral zeros for  ergodic flows }
\author{\bf Valery V. Ryzhikov}
\date{}
\usepackage[T1]{fontenc} 
\usepackage[cp1251]{inputenc}  
\usepackage[russian]{babel}
\usepackage{graphicx}
\graphicspath{{pictures/}}
\DeclareGraphicsExtensions{.pdf,.png,.jpg}
\begin{document}

\maketitle
\Large
   \it  Let a flow $T_t$ preserve an ergodic probability measure $\mu$, $\int f\,d\mu=0$, and $\mu(A)>0$. Then for almost all $x\in A$, for which $f(x)\neq 0$, there is a sequence ${t_k}\to \infty$ such that $T_{t_k}x\in A$ and $\int_0^{t_k} f(T_sx)ds=0$. \rm

\section{ Теоремы Крыгина-Аткинсона и Шнейберга} 
 Эргодическая теория начинается с теоремы Пуанкаре о возвращении.
\it Пусть $T$ -- сохраняющее вероятностную меру $\mu$ обратимое преобразование (автоморфизм)  пространства $(X, \mu)$, тогда для всякого $A\subset X$,
$\mu(A)>0$, почти все точки  $x\in A$ возвращаются в $A$: для некоторого $n>0$ выполняется  $ T^nx\in A.$\rm

Теорема  имеет множество приложений.
Например, существование периодических траекторий в рациональном многоугольнике 
является ее красивым следствием \cite{VGS}. 
Другое начало метрической (от слова мера)  теории динамических систем --  эргодическая теорема Биркгофа (см. \cite{KSF}, \cite{K}), которая, в частности, утверждает, что 
\it для эргодического автоморфизма $T$ и функции $f\in L_1(\mu)$
для п.в. $x$ выполняется
$$\frac 1 n \sum_{i=0}^{n-1}f(T^ix) \ \to\ \int_Xf\,d\mu.$$\rm 
Аналогичным образом  теорема Биркгофа формулируется для сохраняющих меру 
эргодических потоков: \it для почти всех $x\in X$
$$\int_{0}^{t}f(T_sx)\,ds \  = t\int_Xf\,d\mu \ +\ o(t), \ t\to\infty.$$ \rm
Кроме того для почти всякого $x$ найдется последовательность $t_k\to\infty$ такая, что 
$$\int_{0}^{t_k}f(T_sx)\,ds \  = t_k\int_Xf\,d\mu.$$ \rm
 Этот факт следует из результата Шнейберга (см. ниже),  являющегося  непрерывным аналогом дискретной теоремы Крыгина-Актинсона. Сформулируем эти результаты.

\vspace{3mm} 
Пусть $T$ -- автоморфизм вероятностного 
 пространства $(X, \mu)$, $f:X\to\R$ -- $\mu$-измеримая функция.  
\it
Цилиндрический каскад \rm $C: X\times \R\to X\times \R$ определен формулой
$$ C(x,z)=\left(Tx, z+f(x)\right),$$
он сохраняет меру $\bar \mu =\mu\times m$,  где $m$ -- мера Лебега на $\R$.
 
Пусть $T_t$ -- поток, сохраняющий меру $\mu$.  \it Цилиндрический поток  
\rm $C_t: X\times \R\to X\times \R$ задается аналогично:
$$ C_t(x,r)=\left(T_tx, r+ \Phi(t,x)\right), \ \
\Phi(t,x):=\int_{0}^{t}f(T_sx)\,ds.$$ 
Он также сохраняет меру  $\bar \mu =\mu\times m$.

Крыгин \cite{Kr} и Аткинсон \cite{At}  доказали, что цилиндрический каскад с эргодическим $T$ и функцией $f:X\to\R$ с нулевым средним  является рекуррентным.  Это означает, что для всякого множества  
$B\subset X\times \R$ положительной меры для почти всех точек множества $B$ верно, что  они под действием степеней цилиндрического каскада бесконечное число раз  возвращаются в $B$.  Из отмеченной  рекуррентности вытекает, что для всякого $\eps>0$ и почти всех $x\in X$ найдется последовательность $n_k\to\infty$, для которой
  $$\left|\sum_{i=0}^{n_k-1}f(T^ix)\right|<\eps.$$
 Если же функция $f$ принимает целые значения и $\eps <1$, то   
$$\sum_{i=0}^{n_k-1}f(T^ix)=0.$$ 
   Крыгин сформулировал свой результат в случае,
когда $T$ -- эргодический поворот, но в доказательстве использовал только свойство эргодичности  автоморфизма  $T$.

В случае, когда $X$  является метрическим пространством и всякое открытое множество имеет положительную меру,  из отмеченной рекуррентности
 вытекает для почти всех  $x\in X$ существование   последовательности 
$n_k\to\infty$, зависящей от $x$, для которой 
$$T^{n_k}x\to x,  \ \ S(x,n_k):=\sum_{i=0}^{n_k-1}f(T^ix)\to 0.$$

 Теорема 2 из заметки Шнейберга \cite{Sh}  дает дополнительную информацию об этой сходимости: найдутся 
 последовательности натуральных чисел $n_k, p_k\to\infty$  такие, что для всех $k$ выполнено 
$$ S(x,n_k) \leq 0, \ \  S(x, p_k) \geq  0, \ \ \ S(x,n_k)\to 0, \ \  S(x, p_k) \to 0.$$ \rm
Интуитивно понятно, что для типичой функции  $f$ с нулевым средним   суммы  Биркгофа не обращаются в   0 и  сходимость к 0   немонотонна. На самом деле немонотонный характер сходимости имеет место  в  общем 
случае, см. \cite{MRW} (на этот препринт обратил наше внимание Т. Адамс).

Мы будем интересоваться поведением интегралов вдоль траекторий потока. В работе  \cite{Sh}
доказано, что \it для эргодического потока $T_t: X\to X$  и $f:X\to\R$ с нулевым средним  для почти всех $x\in X$  найдется 
последовательность $t_k\to\infty$ такая, что 
$ \int_0^{t_k} f(T_sx)\, ds =0.$ \rm

Для  фиксированного потока $T_t$ и заданной функции  $f$  такие точки  $T_{t_k}x$  будем назвать \it интегральными нулями \rm точки $x$.

\section{ Возвращаемость интегральных нулей для эргодического потока }
Если пространство $X$ дополнительно  оснащено метрикой,  возникает 
   задача о поведении интегральных $x$-нулей в окрестности точки $x$.

\vspace{3mm}
\bf Теорема  (Денисова \cite{De}). \it Пусть  $T_t$ -- эргодический поток на компактном  метрическом пространстве $X$  с конечной мерой Каратеодори $\mu$, а  функция  $f:X\to\R$  непрерывна и имеет  нулевое среднее значение. Тогда для почти всех $x$ при  $f(x)\neq 0$  найдется  последовательность $t_k\to\infty$, для которой  $\int_0^{t_k} f(T_s x)\,ds=0, \ \ T_{t_k}x\to x.$\rm

\vspace{3mm}
Эта теорема и доклад В.В. Козлова   "Возвращаемость интегралов условно периодических функций и эргодическая теория" \    на семинаре по теории функций действительного переменного (мехмат МГУ, 2023)   послужили поводом для нахождения   общего  утверждения, свободного от метрики.

\vspace{3mm}
\bf Теорема 1. \it Пусть  $T_t$ -- сохраняющий меру измеримый   эргодический поток на вероятностном пространстве $(X,\mu)$,    функция $f:X\to\R$  имеет    нулевое  среднее и для $A\subset X$ выполнено  $\mu(A)>0$. Тогда для почти всех $x\in A$ при  $f(x)\neq 0$   найдется  последовательность $t_k\to\infty$, для которой $$\Phi(t_k,x):=\int_0^{t_k} f(T_s x)\,ds=0, \ \ T_{t_k}x\in A. $$\rm

\vspace{3mm}
\bf Следствие.
\it  Пусть  $T_t$ -- эргодический относительно меры $\mu$ 
поток на   метрическом пространстве $X$, мера $\mu$ каждого открытого множества положительна.  Пусть  функция $f:X\to\R$  имеет    нулевое  среднее, тогда для почти всех $x$ при  $f(x)\neq 0$   найдется  последовательность $t_k\to\infty$, для которой $\Phi(t_k,x)=0$ и $T_{t_k}x\to x$.\rm

\vspace{3mm}
 Измеримый поток, сохраняющий меру,  изоморфен специальному потоку $T_t$.
Мы рассматриваем невырожденный случай, когда орбиты всех точек имеют меру 0. 
(В случае потока, являющегося вращением  окружности, доказательство  теоремы тривиально, $\Phi(t,x)$ является непрерывной периодической функцией.) 
Фазовое пространство  $X$  специального потока $T_t$  является частью плоскости под графиком интегрируемой неотрицательной функции $r: [0,1] \to \R^+$. 
Напомним, что точка $x=(a,b)$ ($a\in [0,1]$, $b\in \R^+$)  фазового пространства движется вертикально вверх с постоянной скоростью: $T_t(a,b)=(a, b+t)$. Верхняя граница фазового пространства склеена с нижей таким образом, что точки $(a, r(a))$ отождествляется с точками $(S(a),0)$, где $S$ -- заданный автоморфизм пространства  $([0,1], m)$, $m$  -- мера Лебега на $[0,1]$. В силу теоремы Рудольфа можно  считать, что функция $r$ принимает только два значения, а фазовое пространство $X$  является объединением двух прямоугольников   (см. \cite{KSF}).
Так как всякий измеримый поток  имеет специальное представление, теорема 1 эквивалентна следующему утверждению. 

\vspace{3mm}
\bf Теорема 2. \it Пусть  $T_t$ -- специальный эргодический поток на $(X,\mu)$ и $\mu(A)>0$.  Для   функции $f:X\to\R$  с   нулевым  средним   для почти всех $x\in A\cap \{x':f(x')\neq 0\}$   найдется  последовательность $t_k\to\infty$, для которой  $\Phi(t_k,x)=0, \ \ T_{t_k}x\in A. $ \rm

\vspace{3mm}
Доказательство.  Мы покажем, что при выполнении условий теоремы  для всякого множества $A\subset  \{x':f(x')\neq 0\}$ положительной меры   и всякого натурального $N$   найдется точка $x$ и число $t_1\geq 0$, $t>N$, для которых 
 $$\Phi(t_1,x)= \Phi({t_1+t}, x), \ \ T_{t_1}x,\ T_{t_1+t}x\in A. $$  
Тем самым утверждение теоремы будет доказано. 

  Так как  $f$ интегрируема на $X$,    из теоремы Фубини следует, что функции  $f(T_sx)$ как функции от $s$  для почти всех $x$  интегрируемы на вертикальных отрезках, из которых составлено фазовое пространство потока.
Поэтому для почти всех $x$ неопределенный интеграл  
$$\Phi(t,x)= \int_0^{t} f(T_sx)\, ds$$ является  абсолютно непрерывной функцией
на отрезках. В этом случае  для почти всех  $t$  существует предел
$$ \lim_{\Delta \to 0}\frac  {\int_t^{t+\Delta}f(T_sx)\, ds}{\Delta}=f(T_tx).$$  Тогда для почти  всех 
  $x$  имеем
$$ \lim_{ t\to 0}\frac  {\int_0^{t}f(T_sx)\, ds}{t}=f(x). $$
  И.В. Подвигин сообщил автору, что утверждение о существовании п.в. такого предела для измеримых потоков известно как локальная теорема Винера, см. \cite{K}.

Для всякого $\delta>0$ при условии, что $f$ ненулевая функция,   найдутся  $c> 0$ и  множество  $A'\subset A$ положительной меры, на котором функция $f$ принимает значения из отрезка $[c, (1+\delta) c]$. Пусть    $0<\delta< 10^{-3}$. Так как интегральные нули потока не изменятся, если вместо $f$ рассмотреть функцию  $f/c$, можно считать, что   $c=1$. 
В силу сказанного выше для почти всех $x\in A'$ найдется   число  $b(x)>0$ такое, что  для всех $t$ при $0< t\leq b(x)$  выполнено 
$$1-\delta <\frac 1 t\int_0^{t}|f(T_sx)| ds < 1+2\delta. $$
При   $b>0$ определим множество $A_b$  как  пересечение
множеств 
$$\left\{x\in A'\,:\, 1-\delta < \frac 1 t\int_0^{t}|f(T_sx)|\, ds< 1+2\delta ,\  0< t\leq b\right\},\eqno(1)$$ 
$$\left\{x\in A'\,:\, 1-\delta < \frac 1 t\int_0^{t}f(T_sx)\, ds< 1+2\delta ,\  0< t\leq b\right\}\eqno(2)$$
и  $$\left\{x\in A'\,:\, 1-\delta < \frac 1 t\int_0^{t}\chi_{A'}(T_sx)\, ds,\  0< t\leq b\right\}.\eqno(3)$$  
Так как для почти всех $x\in A'$ выполняется  
$$\lim_{t\to 0}\frac 1 t \int_0^{t}|f(T_sx)|ds\ =\lim_{t\to 0}\frac 1 t \int_0^{t}f(T_sx)ds=\ f(x)$$ (и это также верно для $f=\chi_{A'}$), получим, что  для  некоторого $b>0$ (фиксируем его)  мера множества $A_b$ положительна.

В силу рекуррентности цилиндрического потока, 
построенного по $T_t$ и функции $f$ (о рекуррентности см. пояснения после формулировки  теоремы 4), для  всякого  $N$ мы найдем  $x\in A_b$ и $t>N$
 такие, что 
$$  T_tx\in A_b, \ \ -\delta b<\Phi(t, x)<0.$$
 Функция $\Phi(t, x)=\int_0^{t} f(T_s x)\,ds$  непрерывна по $t$ и выполнено неравенство
$$\Phi(t+2\delta b, x)  - \Phi(t, x) > \delta b/2,$$
так как $T_tx\in A_b$. Отсюда вытекает, что   $\Phi(t', x)=0$
для некоторого $t'>t$ при условии  $t'<t+2\delta b$.   Фиксируем указанное  $t'$.
Покажем, что найдутся такие $s, s'$, $0\leq s,s' <b$,  что 
$$ \Phi(s, x)=\Phi(s', T_{t'}x),\ \ T_{s}x, T_{t'+s'}x \in A_b.\eqno (4)$$
Тогда  получим $\Phi(t'+s'-s, T_{s}x)=0$, и  тем самым теорема будет доказана.

Рассмотрим  множества 
$$ B=\{s\in [0,b]\,:\, T_sx\in A_b\},\ \ B'=\{s'\in [0,b]\,:\, T_{t'+s'}x\in A_b\}.$$ 
Обозначим через $m$ стандартную меру Лебега на прямой. Из (3) получаем 
$$\int_Bf(T_sx)\, ds\ > \ m(B)> (1-\delta)b,$$
так как на множестве $B$ (и на $B'$) функция $f$ превосходит 1.
Обозначим  $D=[0,b]\setminus B$ и $D'=[0,b]\setminus B'$.
Из (1) и равенства 
$$\int_{[0,b]}|f(T_sx)|\, ds =\int_Bf(T_sx)\, ds +\int_D|f(T_sx)|\, ds$$
получаем 
$$\int_D|f(T_sx)|\, ds\ < \ 3\delta b\eqno (5)$$
и аналогичное неравенство для множества $D'$.

\vspace{3mm}
\bf Лемма 1. \it Пусть $f\in L_1[0,b]$,\ $b>0$, \ $F(t)=\int_0^tf(s)\,ds$.  Тогда для измеримого множества $D\subset [0,b]$  выполнено
неравенство $ m(F D)\leq \int_D|f|\,dm.$\rm

\vspace{3mm}
Рассмотрим функции  
$F_1(s)= \Phi(s, x),$ $F_2(s')= \Phi(s', T_{t'}x).$ 
Из  (2) имеем 
$(1-\delta) b<F_1(b), F_2(b) < (1+2\delta) b$.  Из леммы 1 и (5) вытекает $m(F_1 D), \  m(F_2 D')\, < 3\delta b.$ При  малом 
$\delta$ множества  $F_1 B, \  F_2 B'$  по мере $m$ мало отличаются от отрезка $[0,b]$ и их пересечение непусто, так как оно имеет меру, близкую к $b$. 
Таким образом, найдутся $s\in B$, $s'\in B'$ такие, что $F_1(s)=F_2(s')$.
Но это влечет за собой (4) и тем самым утверждение   теоремы  2.

\vspace{3mm}
Доказательство леммы 1. Для  всякого $d>0$  найдется  открытое множество $G$, 
накрывающее  $D$, причем   $m(G)< m(D) +d$. Покажем, что 
$$ m(F G)\leq \int_G|f|\,dm.\eqno (6)$$
Множество $G$ является объединением непересекающихся интервалов.
Пусть $(t_1, t_2)$ -- один из них.   
Для всякого $t\in (t_1, t_2)$ выполнено
$$\int_{t_1}^t f dm\leq  \int_{t_1}^t f_{+} dm, $$
где $f_{+}$ определена как   $f_{+}(t)=f(t)$ при $f(t)\geq 0$,  $f_{+}(t)=0$
при  $f(t)< 0$.
Также выполнено  
$$\int_{t_1}^t f dm\geq  \int_{t_1}^t f_{-} dm, $$
где  $f_{-}(t)=f(t)$ при $f(t)< 0$, иначе $f_{-}(t)=0$.
Так как образ $[t_1, t_2]$ при отображении $F$ есть 
$$[\min \{F(t): t_1\leq t\leq t_2\}\,,\, \max \{F(t): t_1\leq t\leq t_2\}],$$   получаем, что его мера не превосходит 
$$\int_{t_1}^{t_2}|f_{-}|\,dm +\int_{t_1}^{t_2}f_{+}\,dm 
=\int_{t_1}^{t_2}|f|\,dm.$$
Таким образом,   (6) установлено.
Имеем  $$ m(F D)\leq \int_G|f|\,dm \ \leq \ \eps(d)\ + \ \int_D|f|\,dm,$$
где $\eps(d)\to 0$ при $d\to 0$ в силу абсолютной непрерывности интеграла. Лемма и теорема 2 доказаны.

\section{ Замечания}
Рекуррентности интегралов вдоль таекторий и изучению множеств нерекуррентности
посвящено большое количество работ, см. \cite{Mo}-\cite{Ko} и ссылки.
Эта тематика непосредственно связана со случайными блужданиями на одномерных группах $\Z,\R$.  Отметим также, что рекуррентность наблюдается в специальных случаях при блужданиях на матричных группах (см. \cite{Li}) и
даже для типичных блужданий на группе всех автоморфизмов фиксированного вероятностного пространства \cite{23}, \cite{23Sb}. 
Можно дополнить теорему Крыгина-Аткинсона следующим утверждением.

\vspace{3mm} 
 \bf Теорема 3.  \it  Пусть $T$ -- эргодический автоморфизм  стандартного пространства $(X,\mu)$ с $\sigma$-конечной мерой $\mu$.  Если $f:X\to\R$  имеет нулевое среднее, то 
цилиндрический каскад  $C: X\times \R\to X\times \R$,
$ C(x,r)=\left(Tx, r+f(x)\right),$ является рекуррентным. \rm

\vspace{3mm} 
Доказательство.  Рассмотрим индуцированный цилиндрический каскад на множестве $A\times \R$,  $\mu(A)=1$, $A\subset X$. 
Пусть $x\in A$, определим время возврата в $A$ равным  $n(x)$, если  выполнено $S^{n(x)}\in A$, $S^k\notin A,$ при $ 0<k<n(x)$.
Положим 
$$\tilde Sx= S^{n(x)}x, \ \ \ \tilde f (x)=\sum_{i=0}^{n(x)-1}f(S^ix),
\ x\in A.$$ 
Так как  $\tilde S:A\to A$  эргодичен на $A$ и $\int_A  \tilde f (x)\, d\mu =0$,
индуцированный цилиндрический каскад $\tilde C:A\times\R\to A\times \R$ рекуррентный, следовательно, таковым является исходный каскад $C$.

\vspace{2mm}
\bf Теорема  4.  \it  Пусть $S_t$ -- эргодический и рекуррентный поток на стандартном пространстве $(X,\mu)$ с $\sigma$-конечной мерой $\mu$. Если $f:X\to\R$  имеет нулевое среднее, то 
цилиндрический поток  $C_t: X\times \R \to X\times \R$,\
$ C_t(x,r)=\left(S_tx, r+\int_0^t f(S_t x)\, ds\right),$ является рекуррентным. \rm

\vspace{2mm}
Предположим, что существует такое $a>0$, что $S_a$  --  эргодический автоморфизм. Пусть $a=1$.
Рассмотрим функцию $$\tilde f(x):=\int_0^1 f(S_t x)\, ds.$$  Так как $\int_X\tilde  f(x)\, ds=0$, получим, что цилиндрический каскад  $$ \tilde C(x,r)=\left(S_1x, r+\tilde f(x)\right),$$ является рекуррентным в силу теоремы 3.  Но его рекуррентность влечет рекуррентность автоморфизма  $C_1$, так как $C_1=\tilde C$, что вытекает из определений. А из рекуррентности каскада $C_1$ непосредственно  следует рекуррентность цилиндрического потока.
Таким образом, нам  следует  убедиться в том, что эргодичность потока $S_t$ плюс его рекуррентность  влечет за собой эргодичность  автоморфизма $S_a$ для некоторого $a>0$.  

В случае пространства с конечной мерой эргодичность является спектральным свойством и все автоморфизмы  $S_a$, кроме, быть может,  счетного множества,  являются эргодическими.  Это, заметим,  доказывает утверждение теоремы в случае, когда $X$ имеет конечную меру.   В  случае бесконечной меры множества $X$ эргодичность потока  не является спектральным свойством.  Но имеет место следующее утверждение.

\vspace{3mm}
\bf Лемма 2. \it Для эргодического и  рекуррентного  потока $\{S_t\}$, действующего на  пространстве с $\sigma$-конечной мерой,  для почти всех $a\in\R$  автоморфизмы $S_a$ являются эргодическими.
\rm  

\vspace{3mm}
Б. Вейс  сообщил автору, что  лемма 2 вытекает известных фактов, а Дж. Ааронсон  указал  на работу  \cite{Sch} по этой тематике.

Отметим, что  условие интегрируемости функции $f$ не обязательно для рекуррентности цилиндрического каскада.   Из теоремы 1.4  \cite{We} вытекает следующее утверждение.  

\vspace{3mm} 
\it  Пусть $(X,\mu)$ -- вероятностное  пространство и $T$ -- его автоморфизм.  Если $f:X\to\Z$  такова, что 
$$\mu\left( x: \left|\sum_0^{n-1} f(T^ix)\right|>\eps n\right)\to 0, \ n\to\infty,$$  то для почти всех $x$  бесконечное число раз сумма $\sum_0^{n-1} f(T^ix)$ принимает нулевое значение.   \rm

\vspace{10mm}
Автор благодарит Дж. Ааронсона,\ Т. Адамса, \ Б. Вейса,\\ И.В. Подвигина и Ж.-П. Тувено за отклик и  полезные замечания.

\newpage
Recurrence of integral zeros for  ergodic flows

V.V. Ryzhikov

Abstract: Let $T_t$ be a measure-preserving measurable ergodic flow on a probability space $(X,\mu)$, let $f:X\to\R$ have zero mean, and let $A\subset X$ satisfy $\mu(A)>0$. Then for almost all $x\in A$ with $f(x)\neq 0$ there exists a sequence $t_k\to\infty$ such that   $ T_{t_k}x\in A $ and 
$\int_0^{t_k} f(T_s x)\,ds=0.$

Key words: cylindrical flow, ergodicity, recurrence, integrals along trajectories.

\end{document}